\documentclass[reqno,11pt,twoside,tbtags]{amsart}

\usepackage[a4paper,centering,twoside]{geometry}

\numberwithin{equation}{section}

\usepackage{amsmath}
\usepackage{amssymb}
\usepackage{amsfonts}
\usepackage{amsthm}
\usepackage{stmaryrd}

\title{Systoles of hyperbolic manifolds}
\author{Mikhail V.~Belolipetsky} 
\address{Department of Mathematical Sciences, Durham~University, South~Road, Durham, DH1~3LE, England, and Institute of Mathematics, Koptyuga 4, 630090~Novosibirsk, Russia}
\email{mikhail.belolipetsky@durham.ac.uk}
\author{Scott A.~Thomson}
\address{Department of Mathematical Sciences, Durham~University, South~Road, Durham, DH1~3LE, England}
\email{s.a.thomson@durham.ac.uk}
\date{March 2011}

\theoremstyle{plain}
\newtheorem{theorem}{Theorem}[section]
\newtheorem{lemma}[theorem]{Lemma}
\newtheorem{proposition}[theorem]{Proposition}

\theoremstyle{definition}

\theoremstyle{remark}
\newtheorem{remark}[theorem]{Remark}

\renewcommand{\geq}{\geqslant} 

\renewcommand{\leq}{\leqslant} 

\DeclareMathOperator{\PO}{PO}
\DeclareMathOperator{\Ort}{O}
\DeclareMathOperator{\Gal}{Gal}
\DeclareMathOperator{\Sp}{Sp}
\DeclareMathOperator{\PU}{PU}
\newcommand{\R}{\mathbb{R}} 
\newcommand{\Q}{\mathbb{Q}} 
\newcommand{\Z}{\mathbb{Z}} 
\newcommand{\Int}{\mathcal{O}}  
\newcommand{\Hyp}{\mathcal{H}} 
\newcommand{\p}{\mathcal{P}} 
\DeclareMathOperator{\Syst}{Syst_{1}} 
\DeclareMathOperator{\Vol}{Vol} 
\DeclareMathOperator{\Isom}{Isom} 
\newcommand{\isom}{\Isom(\Hyp^{n})} 
\newcommand{\embeds}{\hookrightarrow} 

\begin{document}

\begin{abstract}
We show that for every $n\geq 2$ and any $\epsilon>0$ there exists a compact hyperbolic $n$-manifold with a closed geodesic of length less than $\epsilon$. When $\epsilon$ is sufficiently small these manifolds are non-arithmetic, and they are obtained by a generalised inbreeding construction which was first suggested by Agol for $n=4$. We also show that for $n\geq 3$ the volumes of these manifolds grow at least as $1/\epsilon^{n-2}$ when $\epsilon\to 0$.
\end{abstract}

\maketitle

\section{Introduction}

Let $\Hyp^n$ denote the hyperbolic $n$-space. By a compact hyperbolic $n$-manifold we mean a quotient space $M = \Gamma\backslash\Hyp^n$ where $\Gamma$ is a cocompact torsion-free discrete subgroup of $\isom$, the group of isometries of $\Hyp^n$. The \emph{systole} of a compact Riemannian manifold $M$, denoted by $\Syst(M)$, is the length of a closed geodesic of the shortest length in $M$. We refer to a recent monograph by M.~Katz \cite{Katz:systoles} for more information about systoles and systolic geometry.

It is well-known that for any $\epsilon>0$ there exist $2$-dimensional compact hyperbolic manifolds having a systole of length less than $\epsilon$, and examples of such manifolds of any genus $g\geq 2$ can be easily constructed using Teichm\"uller theory. A similar result for $n = 3$ can be achieved using Thurston's hyperbolic Dehn surgery theorem \cite[Theorem 5.8.2]{Thurston:GT3M}. For a long time the existence of compact hyperbolic manifolds with arbitrarily short systoles in higher dimensions was an open problem. In a recent paper \cite{Agol:systoles}, Agol suggested a very interesting construction which solves the problem for $n = 4$. His paper was a starting point for our work.

Our main result is the following:

\begin{theorem}\label{thm:existence}
For every $n\geq 2$ and any $\epsilon>0$, there exist compact $n$-dimensional hyperbolic manifolds $M$ with $\Syst(M)<\epsilon$.
\end{theorem}

The manifolds $M$ are obtained by a variant of an \emph{inbreeding} construction which was first suggested by Agol for $n=4$. We simplify the critical step in the argument of \cite{Agol:systoles} which makes extensive use of geometrical finiteness and related properties, and hence limits his construction to some special examples. A principal ingredient in our proof is a lemma of Margulis and Vinberg \cite{MV:lineargroups} which we generalise to cocompact discrete subgroups of $\isom$. The proof of this generalised Margulis-Vinberg lemma is the main technical part of the proof of the theorem.

\medskip

Systolic geometry studies relations between systole length and volume captured by isosystolic inequalities (cf.~\cite{Katz:systoles}). Our second result provides an inequality of this type for the manifolds from the proof of Theorem~\ref{thm:existence}. To put it into context recall that for $n\geq 4$ (in contrast with $n=2$ and $3$) there exist only finitely many non-isometric hyperbolic $n$-manifolds of bounded volume \cite{Wang}. Hence when $\epsilon\to 0$ for these dimensions we will necessarily have $\Vol(M)\to\infty$. It is natural to ask how fast the volume grows, and the following theorem gives the answer to this question for our manifolds.

\begin{theorem}\label{thm:volumegrowth}
For every $n\geq 3$ there exists a positive constant $C_n$ (which depends only on $n$), such that the systole length and volume of the manifolds obtained in the proof of Theorem~\ref{thm:existence} satisfy
	\[ \Vol(M) \geq \frac{C_n}{\Syst(M)^{n-2}}. \]
\end{theorem}
We can compare this inequality with the well known Gromov isosystolic inequality (\cite[Theorem 12.2.2]{Katz:systoles}), which implies
\[\Vol(M) \geq A_n \Syst(M)^n,\]
where $A_n$ is a positive constant which depends only on $n$. When the systole of $M$ is sufficiently large Gromov's inequality gives a better bound for the volume; however, when the systole is small the inequality of Theorem~\ref{thm:volumegrowth} becomes stronger.

The proof of this theorem uses important recent work of Bridgeman and Kahn on orthospectra and volumes of hyperbolic $n$-manifolds \cite{BK:Orthospectra}.
In fact, we can show that it is possible to achieve that $\Vol(M)$ grows exactly like a polynomial in $1/\Syst(M)$ (see the discussion after the proof of Theorem~\ref{thm:volumegrowth} and Proposition~\ref{prop:volumegrowth}). Therefore, Theorem~\ref{thm:volumegrowth} captures the growth rate of the volume in our construction. 

The paper is organised as follows: In Sect.~\ref{sec:2}, we prove Theorem~\ref{thm:existence} modulo the generalised Margulis-Vinberg lemma. The proof of the lemma is given in Sect.~\ref{sec:3}. The next section is dedicated to Theorem~\ref{thm:volumegrowth}. We end with remarks regarding arithmeticity, commensurability and related questions in Sect.~\ref{sec:5}.

\medskip

{\bf\noindent Acknowledgements.} We thank Ian Agol and Misha Kapovich for helpful discussions and suggestions. We are also grateful to the referee for careful proofreading of the manuscript.
M.~Belolipetsky is partially supported by EPSRC grant EP/F022662/1.

\section{Proof of Theorem 1.1}\label{sec:2}

Let $\epsilon>0$. We first obtain a hyperbolic $n$-manifold $M$ which contains a geodesic segment of length at most $\epsilon/2$ that is orthogonal to two hypersurfaces.

Fix a totally real number field $K\subset\R$, and let $f$ be a non-singular $(n+1)$-ary quadratic form of signature $(n,1)$ defined over $K$. We assume that for Galois embeddings $\sigma\colon K\embeds\R$ which are different from the original embedding, $f^\sigma$ is positive definite. It is well known that $\PO_{f}(\R)$ is isomorphic to the group of isometries $\isom$ of the hyperbolic $n$-space, and that $\PO_{f}(\Int_{K})$ (${=} \Ort_{f}(\Int_{K})/\{+1,-1\}$) is an arithmetic lattice in $\PO_{f}(\R)$, where $\Int_{K}$ denotes the ring of integers of $K$.
The lattices obtained this way (and subgroups of $\PO_{f}(\R)$ which are commensurable with them) are called \emph{arithmetic subgroups of the simplest type}. From now on let us assume that the degree of the field $K$ is at least $2$; i.e., $K\neq\Q$. Then by Godement's compactness criterion $\PO_{f}(\Int_{K})$ is cocompact. Now, by Selberg's Lemma, one can find a torsion-free subgroup $\Gamma < \PO_{f}(\Int_{K})$, of finite index. Thus $\Gamma\backslash\Hyp^{n}$ is a compact hyperbolic $n$-manifold. 
(We refer to \cite{MR:book} for more details about arithmetic subgroups and their properties.)

Consider the vector model of hyperbolic $n$-space associated with $f$. More precisely, we can define an inner product on $\R^{n+1}$ by
\begin{equation}\label{eq:product}
(u,v) = (u,v)_f = \tfrac12\bigl(f(u+v)-f(u)-f(v)\bigr).
\end{equation}
After scaling we can assume that all coefficients of $f$ are in $\Int_{K}$, and, moreover, for any $u,v\in\Int_{K}^{n+1}$, we have $(u,v)_{f}\in\Int_{K}$.
This does not affect the generality of the construction and will be our standing assumption.
The form $f$ has signature $(n,1)$ which implies that $(\cdot,\cdot)_f$ is equivalent to a standard Lorentzian inner product on $\R^{n+1}$, and therefore we can repeat the well-known construction of the vector model of $\Hyp^n$ using the inner product $(\cdot,\cdot)_f$ (cf.~ \cite[Ch.~3]{Ratcliffe:Foundations}).

Let us choose a vector $e_{0}\in K^{n+1}$ with $f(e_0)>0$. The intersection $H_{0}=\langle e_0 \rangle^{\perp_f}\cap\Hyp^{n}$ in the ambient space $\R^{n+1}$ is a hyperplane in $\Hyp^{n}$ and, moreover, $\Gamma_0 = \mathrm{Isom}(H_0)\cap\Gamma$ is a cocompact discrete subgroup of $\mathrm{Isom}(H_0)$ (where we embed $\mathrm{Isom}(H_0)\hookrightarrow\isom$ in the natural way). The latter holds because $\Gamma$ is defined over $K$ and $H_0$ is a $K$-rational hyperplane. 
Now we would like to find another $K$-rational hyperplane $H_1$ which is $\epsilon/2$-close to $H_0$ but disjoint from it. This can be achieved as follows.
A well-known result of Margulis asserts that since the group $\Gamma$ is arithmetic, its commensurator is dense in $\PO_{f}(\R)$ \cite[Theorem~1, p.~2]{Margulis:book}.
Indeed, this commensurator contains $\PO_{f}(K)$, so we can find $\gamma\in \PO_{f}(K)$ with the property that $H_{1}$~--- where $H_{1}=\langle e_{1} \rangle^{\perp_f} \cap \Hyp^{n}$ and $e_{1}=\gamma(e_{0})$~--- is disjoint, but at most distance $\epsilon/2$ away, from $H_{0}$. The hyperplane $H_1$ is also $K$-rational so $\Gamma_1 = \mathrm{Isom}(H_1)\cap\Gamma$ is a cocompact discrete subgroup of $\mathrm{Isom}(H_1)$, and both $\Gamma_0\backslash H_0$ and $\Gamma_1\backslash H_1$ are immersed into $\Gamma\backslash \Hyp^{n}$ as totally geodesic hypersurfaces.

The generalised Margulis-Vinberg lemma (below) states that one can find a finite-index subgroup $\Gamma'<\PO_{f}(\Int_{K})$ with the property that for every $h\in\Gamma'$,
\[ \text{either } \quad h(H_{0})=H_{0} \quad\text{ or }\quad h(H_{0}) \cap (H_{0}\cup H_{1})=\emptyset. \]
Reversing the roles of $H_0$ and $H_1$ we can apply the lemma again to obtain a subgroup $\Gamma''$ such that for every $h\in\Gamma''$,
either $h(H_{1})=H_{1}$ or $h(H_{1}) \cap (H_{0}\cup H_{1})=\emptyset$.
 
Let $\Lambda = \Gamma\cap\Gamma'\cap\Gamma''$. Then the natural projections of $H_{0}$ and $H_{1}$ in the quotient $\Lambda\backslash\Hyp^{n}$ will be properly embedded and will not intersect. Thus the manifold $L = \Lambda\backslash\Hyp^{n}$ contains properly embedded totally geodesic hypersurfaces $\Lambda_0\backslash H_0$ and $\Lambda_1\backslash H_1$  (where $\Lambda_i = \mathrm{Isom}(H_i)\cap\Lambda$, $i = 0,1$) which are $\epsilon/2$-close.
Let $g$ be a geodesic segment orthogonal to both of them, so as to fulfil our aim stated at the beginning of the proof.

To complete the proof, we `cut' $L$ along the hypersurfaces $\Lambda_0\backslash H_{0}$ and $\Lambda_1\backslash H_{1}$. Retaining the connected component containing $g$ (if the cutting separates the manifold), we have a manifold $M'$ with a totally geodesic boundary and $g$ orthogonal to this boundary. Taking the double of $M'$ results in a closed $n$-manifold $M$, and the segment $g$ becomes a closed geodesic of $M$ of length at most $\epsilon$.
This concludes the construction. \hfill\qed

\section{Generalised Margulis-Vinberg lemma}\label{sec:3}

An earlier form of the result in this section appeared in a paper by Margulis and Vinberg \cite{MV:lineargroups}, and was also used by Kapovich, Potyagailo and Vinberg in a paper on non-coherence of lattices \cite{KPV:noncoherence}.
The generalised version here considers the case where a quadratic form is defined over a number field $K/\Q$ rather than being purely rational. This generalisation is necessary for dealing with cocompact lattices, and leads to a proof that uses an essentially different argument from that for $K=\Q$.

As before, let $K\subset\R$ be a totally real algebraic number field of degree $d$, and $f$ a quadratic form over $K$ of signature $(n,1)$ and such that $f^{\sigma_{j}}$ is positive definite for all Galois embeddings $\sigma_j\colon K\embeds\R$, $j = 2, \dots, d$ which are different from the original embedding (denoted by $\sigma_1$). Thus the group $\PO_{f}(\Int_{K})$ is an arithmetic lattice in $\PO_{f}(\R)$.

Let $H_{0}, \ldots, H_{k}$ be pairwise disjoint hyperplanes in $\Hyp^{n}$ defined by $H_{i}=\Hyp^{n}\cap \langle e_{i}\rangle^{\perp_f}$,  $e_{i}\in K^{n+1}$, $i=0,\ldots,k$.

\begin{lemma}\label{lemma}
	There exists a finite index subgroup $\Gamma' < \PO_{f}(\Int_{K})$ such that for every $ h\in\Gamma',$
\[ \text{either } \quad h(H_{0})=H_{0} \quad\text{ or }\quad h(H_{0}) \cap (H_{0}\cup \cdots \cup H_{k})=\emptyset. \]
\end{lemma}

\begin{proof}
	Write $\Gamma=\PO_{f}(\Int_{K})$.
Since multiplying the $e_{i}$ by the denominators of their entries preserves their orthogonality with the $H_{i}$, we can assume the $e_{i}$ have entries in $\Int_{K}$.
Thus, if $h\in \Gamma$, then $(h(e_{0}),e_{i})\in\Int_{K}$.
(Throughout this proof, inner products and orthogonal complements are understood to be with respect to $f$ as in \eqref{eq:product}.)

\par Assume that $k=1$.
\par Let $\p$ be the principal ideal $\p=(\beta)\subset\Int_{K}$, where $\beta=2C(e_{0},e_{1})$ and $C>1$ is an integer to be determined later.
The congruence subgroup $\Gamma(\p)<\Gamma$ contains (by definition) precisely those elements $h\in\Gamma$ with the property that
\[ h \equiv \mathbf{id} \pmod \p. \]
(Here we use the fact that $\PO_{f}$ can be identified with the matrix group $\Ort'_{f}$, the subgroup of the orthogonal group $\Ort_{f}$ which preserves the upper half-space.)
Hence, given $h\in\Gamma(\p)$, we find that
\begin{equation}\label{eqn:Gamma-action}
	\bigl( h(e_0), e_{1} \bigr) = \bigl(e_{0},e_{1}\bigr) + \alpha\beta
\end{equation}
for some $\alpha\in\Int_{K}$ (where $\alpha$ depends on $h$).
We wish to show that for every $h\in \Gamma_{1}=\Gamma(\p)$ we have $h(H_{0})\cap H_{1}=\emptyset$.

\par To be able to examine the intersections of the hyperplanes we use the following inequality \cite[Theorem~3.2.6]{Ratcliffe:Foundations}: two hyperplanes in $\Hyp^{n}$, defined as above by normal vectors $v_{0}$ and $v_{1}$, intersect transversally if and only if
\[ \bigl| (v_{0},v_{1}) \bigr| < \sqrt{ (v_{0},v_{0})(v_{1},v_{1}) }. \]
So we find that the hyperplanes are disjoint or coincide completely if
\begin{equation}\label{eqn:disjointnesscondition}
\bigl| (v_{0},v_{1}) \bigr| \geq \sqrt{ (v_{0},v_{0})(v_{1},v_{1}) }.
\end{equation}

\par If $\alpha=0$ in \eqref{eqn:Gamma-action} then
\begin{equation}\label{eqn:alphaiszero}
  \bigl| \bigl( h(e_{0}),e_{1} \bigr) \bigr| =  \bigl| \bigl( e_{0},e_{1} \bigr) \bigr| \geq \sqrt{\vphantom{\bigl(}(e_{0},e_{0})(e_{1},e_{1})} = \sqrt{\bigl(h(e_{0}),h(e_{0})\bigr)(e_{1},e_{1})}
\end{equation}
(where the inequality follows from \eqref{eqn:disjointnesscondition} and the initial condition $H_{0}\cap H_{1}=\emptyset$), and hence the hyperplanes $h(H_{0})$ and $H_{1}$ are either disjoint or equal.
We will eliminate the possibility of equality later in the proof.

\par If $|\alpha|\geq 1$, then
\begin{align*}
\bigl| \bigl( h(e_{0}),e_{1} \bigr) \bigr| &= \bigl| (e_{0},e_{1}) + \alpha\beta \bigr| \\
&\geq \bigl| |\alpha||\beta| - |(e_{0},e_{1})|\bigr| &&\text{(by the triangle inequality for $|\cdot|$)}\\
	&= \bigl| (e_{0},e_{1}) \bigr|\cdot \bigl| 2C|\alpha| - 1 \bigr| \\
	&> \bigl| (e_{0},e_{1}) \bigr| &&\text{(since $C>1$)}.
\end{align*}

\par The case $|\alpha|<1$ requires more care.
Writing $x^{\sigma}$ for the conjugate of $x\in K$ by $\sigma\in\Gal(K/\Q)$, we have the norm $N(x)$ of $x$ defined as the product $\prod_{\sigma\in\Gal(K/\Q)}x^{\sigma}$ \cite[Ch. I, Sect.~5]{Lang:ANT}; and if $x\in\Int_{K}$ then $N(x)\in\Z$ so that $|N(x)|\geq 1$ for $x\in\Int_{K}\setminus\{0\}$ \cite[Ch.~I, Cor. to Prop.~5]{Lang:ANT}.
Since $\alpha\in\Int_{K}\setminus\{0\}$ and $|\alpha|<1$, the preceeding means that $|\alpha^{\sigma_j}|>1$ for some $j\in\{2,\ldots,d\}$.
For this $j$, we get
  \begin{align*}
 \bigl| \bigl( h(e_{0}),e_{1} \bigr)^{\sigma_j} \bigr| &= \bigl| (e_{0},e_{1})^{\sigma_j} + \alpha^{\sigma_j}\beta^{\sigma_j} \bigr| \\
	&\geq \Bigl| |\alpha^{\sigma_j}||\beta^{\sigma_j}| - |(e_{0},e_{1})^{\sigma_j}| \Bigr| = \bigl|\beta^{\sigma_j} \bigr|  \left| |\alpha^{\sigma_j}| - \frac{1}{2C} \right| \\
	&> \tfrac{1}{2} \bigl| \beta^{\sigma_j}\bigr| = C \bigl| (e_{0},e_{1})^{\sigma_j} \bigr| .
\end{align*}
We have
\begin{equation}\label{eqn:estimate1}
\bigl| \bigl( h(e_{0}),e_{1} \bigr)^{\sigma_j} \bigr| >   C \bigl| (e_{0},e_{1})^{\sigma_j} \bigr|   =   C \bigl| (e_{0},e_{1})^{\sigma_j} \bigr| \frac{\bigl| \bigl( h(e_{0}),e_{1} \bigr)^{\sigma_j}  \bigr|}{\bigl| \bigl( h(e_{0}),e_{1} \bigr)^{\sigma_j}  \bigr|}. 
\end{equation}
Now since $(\cdot,\cdot)_{f^{\sigma_j}}$ is positive definite (for $j\geq 2$), the Cauchy-Schwarz inequality applies and we can use it to bound the denominator from above.
(Note that we have $(u,v)^{\sigma}=(u^{\sigma},v^{\sigma})_{f^\sigma}$.)
Thus \eqref{eqn:estimate1} becomes
\[ \bigl| \bigl( h(e_{0}),e_{1} \bigr)^{\sigma_j} \bigr| > C \frac{\bigl| (e_{0},e_{1})^{\sigma_j}\bigr| \bigl| \bigl( h(e_{0}),e_{1} \bigr)^{\sigma_j} \bigr| }{\sqrt{\bigl(h(e_{0}),h(e_{0})\bigr)^{\sigma_j}(e_{1},e_{1})^{\sigma_j}}} = C \frac{\bigl| (e_{0},e_{1})^{\sigma_j}\bigr| \bigl| \bigl( h(e_{0}),e_{1} \bigr)^{\sigma_j} \bigr| }{\sqrt{\vphantom{\bigl(h\bigr)^{\sigma_j}} (e_{0},e_{0})^{\sigma_j}(e_{1},e_{1})^{\sigma_j}}} . \]
Multiplying each side of this inequality by all the $\bigl| \bigl( h(e_{0}),e_{1} \bigr)^{\sigma_k} \bigr|$ for which $k\neq j$ gives
\begin{equation}\label{eqn:estimate2}
  \bigl| N\bigl( (h(e_{0}),e_{1}) \bigr) \bigr| > C  \underbrace{\vphantom{\int} \bigl|  N\bigl( (h(e_{0}),e_{1}) \bigr)\bigr| }_{(\ast)} \frac{\bigl| (e_{0},e_{1})^{\sigma_j} \bigr|}{\sqrt{(e_{0},e_{0})^{\sigma_j}(e_{1},e_{1})^{\sigma_j}}} .
\end{equation}
We can replace $(\ast)$ by $\bigl| N\bigl( (e_{0},e_{1}) \bigr) \bigr|$, for
\begin{align}
 \bigl| N\bigl( (h(e_{0}),e_{1} )\bigr) \bigr| &= \bigl| N\bigl( (e_{0},e_{1}) +2C\alpha (e_{0},e_{1}) \bigr) \bigr| \notag \\  
  &= \bigl|N\bigl( (e_{0},e_{1}) \bigr)  \bigr|\cdot \bigl| N(1+2C\alpha) \bigr| &&\text{(by multiplicativity of the norm)}  \notag \\
  &\geq \bigl|N\bigl( (e_{0},e_{1}) \bigr)\bigr| &&\text{(since $1+2C\alpha\in\Int_{K}\setminus \{0\}$).}\label{eqn:estimate3}
\end{align}
Writing both norms in \eqref{eqn:estimate2} as products of Galois conjugates, \eqref{eqn:estimate2} and \eqref{eqn:estimate3} give
\begin{align}
\prod_{i=1}^{d} \bigl| \bigl( h(e_{0}),e_{1} \bigr)^{\sigma_i} \bigr| & >C \left( \prod_{i=1}^{d}\bigl| (e_{0},e_{1})^{\sigma_i} \bigr| \right)  \frac{\bigl|(e_{0},e_{1})^{\sigma_j}\bigr|}{\sqrt{(e_{0},e_{0})^{\sigma_j}(e_{1},e_{1})^{\sigma_j}}} , \notag \\
\intertext{so that by rearranging,}
  \bigl| \bigl( h(e_{0}),e_{1} \bigr) \bigr| &> C \frac{\bigl|(e_{0},e_{1})^{\sigma_j}\bigr|}{\sqrt{(e_{0},e_{0})^{\sigma_j}(e_{1},e_{1})^{\sigma_j}}} \left(\prod_{i=2}^{d} \frac{\bigl| (e_{0},e_{1})^{\sigma_i}\bigr|}{\bigl| \bigl( h(e_{0}),e_{1} \bigr)^{\sigma_i} \bigr| } \right) \bigl| (e_{0},e_{1}) \bigr|\notag  \\
  &\geq   C \frac{\bigl|(e_{0},e_{1})^{\sigma_j}\bigr|}{\sqrt{(e_{0},e_{0})^{\sigma_j}(e_{1},e_{1})^{\sigma_j}}} \left(\prod_{i=2}^{d} \frac{\bigl| (e_{0},e_{1})^{\sigma_i}\bigr|}{\sqrt{ (e_{0},e_{0})^{\sigma_i}(e_{1},e_{1})^{\sigma_i} } } \right) \bigl| (e_{0},e_{1}) \bigr| \label{eqn:finalestimate} .
\end{align}

\par At this point by choosing $C$ to be sufficiently large we can ensure that $\bigl|\bigl(h(e_{0},e_{1})\bigr)\bigr|>\bigl|(e_{0},e_{1})\bigr|$.
Notice that since $\alpha$ depends on $h$, so too does $j$; however, we can choose $C$ independently of $h$ by assuming
\begin{equation}\label{eqn:Cchoice} 
  C \geq \prod_{\sigma_{j}\in S} \frac{(e_{0},e_{0})^{\sigma_j} (e_{1},e_{1})^{\sigma_j} }{[(e_{0},e_{1})^{\sigma_j}]^{2} } ,
\end{equation}
where $S\subseteq\{\sigma_{2},\ldots,\sigma_{d}\}$ is the set of all $\sigma_{j}$ for which the corresponding factors in \eqref{eqn:Cchoice} are greater than 1.
\par Thus we get $\bigl| \bigl( h(e_{0}),e_{1} \bigr) \bigr| \geq \bigl| (e_{0},e_{1}) \bigr| \geq \sqrt{(h(e_{0}),h(e_{0}))(e_{1},e_{1})}$ as in the other two cases for $\alpha$.
This means that $h(H_0)$ either coincides with, or does not intersect $H_{1}$.

\par To avoid the possibility of $h(H_{0})$ coinciding with $H_{1}$, we have to ensure that $h(e_{0}) \neq \pm \omega e_{1}$ for some $\omega \in \R_{>0}$.
If it exists, then this $\omega$ would be given by $\omega = \sqrt{ (e_{0},e_{0}) / (e_{1},e_{1}) }$ and there are two possible cases:
\par (a) $\omega \notin K$, whence $h(e_{0})=\pm\omega e_{1}$ is impossible.
\par (b) $\omega \in K$. 
Let $e_{1}'$ be the vector obtained by scaling $\omega e_{1}$ by $\sqrt{(e_0,e_0)(e_1,e_1)}$, so that $e_{1}'=(e_0,e_0)e_1$.
Similarly, define $e_{0}' = \sqrt{(e_0,e_0)(e_1,e_1)}\, e_{0}$.
Thus we have $h(e_0)=\pm\omega e_{1}$ if and only if $h(e_{0}') = \pm e_{1}'$.
We verify that $e_{0}'$ and $e_{1}'$ are in fact in $\Int_{K}^{n+1}$, using the following:
(1) If $x\in\Int_{K}$ and $\sqrt{x}\in K$ then $\sqrt{x}\in\Int_{K}$.
Indeed, if $x^{m}+a_{m-1}x^{m-1} + \cdots + a_{0} =0$ is an equation (with coefficients in $\Z$) giving rise to the algebraic integer $x$ \cite[Ch.~I, Sect. 2]{Lang:ANT}, then $(\sqrt{x})^{2m} + a_{m-1}(\sqrt{x})^{2(m-1)} + \cdots + a_{0}=0$ holds and hence $\sqrt{x}$ is an algebraic integer.
Denote by $\mathbf{A}$ the set of all algebraic integers, and note that $\Int_{K}=\mathbf{A}\cap K$ (cf.\ \cite[Ch.~1, Prop.~5]{Lang:ANT}).
Thus since $\sqrt{x}\in K$ we have $\sqrt{x}\in\Int_{K}$.
(2) If $\sqrt{x/y}\in K$ (for $x,y\in K$), then $\sqrt{\vphantom{/}xy}\in K$, for $\sqrt{\vphantom{/} xy} = y\sqrt{x/y}\in K$.
\par Now for $h(e_{0}')=\pm e_{1}'$ to hold we must have $e_{0}' + v = \pm e_{1}'$ where $v\equiv 0$ modulo $\p$.
If $e_{0}' + e_{1}'$ and $e_{0}'- e_{1}'$ are not congruent to $0$ modulo $\p$, then this coincidence will not occur, and we can ensure this by choosing $C$ sufficiently large.

\par It remains to check that $h(H_0)$ and $H_0$ either coincide or are disjoint.
One can repeat all of the above argument as far as \eqref{eqn:finalestimate}, with $e_0$ in place of $e_1$, and we find that the ideal $\p'=2(e_0,e_0)$ actually suffices in place of $\p$ to ensure that we have $|(h(e_{0}),e_0)|\geq \sqrt{(h(e_0),h(e_0))(e_0,e_0) }$ for every $h\in\Gamma(\p')$.
Denote $\Gamma(\p')$ by $\Gamma_{0}$. 

\par 
If $k\geq 2$, then to separate all hyperplanes we apply the above argument to all other $e_{i}$, ($i=2,\ldots,k$) so that we get $\Gamma_{2},\ldots,\Gamma_{k}$ which are also finite-index subgroups of $\Gamma$.
The group $\Gamma'=\Gamma_{0} \cap \Gamma_{1}\cap \cdots \cap \Gamma_{k}$ will then satisfy the conclusion of the lemma, and is still of finite index in $\Gamma$.
\end{proof}

\begin{remark}\label{rem:strictness}
  If we assume that the hyperplanes $H_{0},\ldots,H_{k}$ are not only disjoint but also do not meet at infinity then the inequality in \eqref{eqn:alphaiszero} becomes strict and the coincidence of $h(H_{0})$ and $H_{i}$ (for $i=1,\ldots,k$) is automatically avoided.
  \end{remark}

\section{Proof of Theorem 1.2}\label{sec:4}

Let $M'$ be a compact hyperbolic $n$-manifold as in Sect.~\ref{sec:2}, for which $M$ is a double. By the construction, $M'$ has a totally geodesic arc of length $\ell = \frac{1}{2}\Syst(M)$ with endpoints in $\partial M'$. This value $\ell$ appears in the orthospectrum of $M'$ as defined in \cite{BK:Orthospectra}. In order to bound the volume of $M'$, and hence of $M$, we can apply the result of Bridgeman and Kahn which relates the volume and orthospectrum of a compact hyperbolic $n$-manifold with non-empty totally geodesic boundary.

Assuming $n\geq 3$, by \cite[Theorem~1]{BK:Orthospectra} we have $\Vol(M')\geq F_n(\ell)$, and by Lemma~9(3) (loc.~cit.), $\lim_{\ell\to 0} \ell^{n-2}F_n(\ell) = K_n$, where $F_n$ is a continuous monotonically decreasing function $\R_{>0} \to \R_{>0}$ and $K_n$ is an explicit positive constant given there. It follows that there exists $K'_n > 0$ such that if $\ell < 1$, then
$$F_n(\ell) \geq \frac{K'_n}{\ell^{n-2}}.$$
Therefore,
$$\Vol(M) = 2\Vol(M') \geq \frac{2^{n-1}K'_n}{\Syst(M)^{n-2}},$$
if $\Syst(M) < 2$.

\par For $\Syst(M)\geq 2$, we can refer to the Kazhdan-Margulis theorem which asserts that there is a constant $A_{n}>0$ such that $\Vol(M)\geq A_{n}$ \cite{KM}.
Hence we can take $C_n = \min(2^{n-1}K'_n, 2^{n-2}A_n)$ and the theorem is proven.
\hfill\qed

\begin{remark}
It was pointed out by the referee that a similar bound holds in a more general setting if the systole of $M$ is hyperbolic.
Indeed, the volume of a Margulis tube about the systole of $M$ can be bounded below by $C/\ell^{n-2}$ in the case where the corresponding element of $\pi_{1}(M)$ has no rotational part.
We refer the reader to Reznikov \cite{Reznikov} for details.
\end{remark}

\par It is interesting to see how close the inequality of Theorem~\ref{thm:volumegrowth} approximates the actual growth of volume in our construction. In Sect.~\ref{sec:2} the desired manifold $M$ is obtained as a double of $M'$, which in turn  appears as a part of the quotient manifold $L = \Lambda\backslash\Hyp^n$. Hence
$$\Vol(M) = 2\Vol(M') \leq 2\Vol(L).$$
The volume of $L$ depends on the index of $\Lambda$ in the initial group $\Gamma$ which can be estimated from the proof of the generalised Margulis-Vinberg lemma in Sect.~\ref{sec:3}. Let us note that by the argument in Sect.~\ref{sec:3} we can always make the index $|\Gamma:\Lambda|$ arbitrarily large, and hence cannot bound the volume of $L$ from above. However, we are rather interested in understanding how small the volume can get when $\epsilon$ in Theorem~\ref{thm:existence} is close to zero and here we can say more.

\begin{proposition}\label{prop:volumegrowth}
For every $n \geq 2$ there exists a sequence of manifolds $\{ M_i \}$ from Theorem~\ref{thm:existence}, such that when $i \to \infty$, $\Syst(M_i) \to 0$ and
	\[ \Vol(M_i) \leq \frac{B_n}{\Syst(M_i)^{\gamma_n}}, \]
where $B_n$ and $\gamma_n$ are positive constants depending on $n$.
\end{proposition}

\begin{proof}
  Consider a sequence of inbred manifolds $M_i$ with $\Syst(M_i) = \epsilon_i \to 0$ when $i\to\infty$. For each $M_i$ we have associated vectors $e_{0}^{(i)}$ and $e_{1}^{(i)}$ defined in Sect.~\ref{sec:2} and an ideal $\p^{(i)}\in\Int_K$ defined in Sect.~\ref{sec:3}. We have $\bigl|N(\p^{(i)})\bigr| = (2 C_i)^{d} \bigl| N\bigl( (e_{0}^{(i)},e_1^{(i)}) \bigr) \bigr|$, where $C_i$ satisfies inequality \eqref{eqn:Cchoice}.
  Recall that in Sect.~\ref{sec:3} the vectors $e_{0}^{(i)}$ and $e_{1}^{(i)}$ are normalised so that they have coordinates in $\Int_K$. This implies that when the angle between $e_0^{(i)}$ and $e_1^{(i)}$ tends to $0$, either $\bigl|N\bigl( (e_{0}^{(i)},e_{1}^{(i)}) \bigr)\bigr|$ or the lower bound for $C_{i}$ will grow, and hence the absolute value of the norm $\bigl|N(\p^{(i)})\bigr| \to \infty$.

\par We give a more concrete example to provide the $\{M_{i}\}$ for the conclusion of the proposition.
Let $K=\Q(\sqrt{5})$ and $f=-\sqrt{5}x_{0}^{2} + x_{1}^{2} + \cdots + x_{n}^{2}$.
The sequence of matrices
\[ A_{i} =  \begin{pmatrix}%
  \frac{i^2 + \sqrt{5}}{i^2 - \sqrt{5}} & 0 & \cdots & 0 & \frac{-2i}{i^2 - \sqrt{5}} \\%
  0 & 1 & & &  0 \\%
  \vdots &  & \ddots & &\vdots \\%
  0 & & & 1 & 0\\%
  \frac{-2i\sqrt{5}}{i^2 - \sqrt{5}} & 0 & \cdots & 0 & \frac{i^2 + \sqrt{5}}{i^{2} - \sqrt{5}}%
\end{pmatrix}\]
can be shown to lie in $\Ort_{f}'(K)$ (the $K$-points of the subgroup of $\Ort_{f}$ that preserves $\Hyp^{n}$), and clearly $A_{i}\to\mathbf{id}$ as $i\to\infty$.
Let $e_{0}=(0,0,\ldots,0,1)$, so that
\[ A_{i}(e_{0}) = \left( \frac{-2i}{i^2 - \sqrt{5}} \, ,\, 0 \, ,\,  \ldots\, ,\,  0 \, , \, \frac{i^2 + \sqrt{5}}{i^2 - \sqrt{5}} \right) \in K^{n+1}. \]
Rescaling $e_{0}$ and $A_{i}(e_0)$, we define
\[ e_{0}^{(i)} = (0,0,\ldots, i^{2}-\sqrt{5}) \quad\text{ and }\quad e_{1}^{(i)} = (-2i, 0 , \ldots, 0 , i^{2}+\sqrt{5}), \]
which give 
\[ \bigl( e_{0}^{(i)} , e_{1}^{(i)} \bigr) = i^{4}-5 \quad\text{ and }\quad \bigl( e_{0}^{(i)}, e_{0}^{(i)} \bigr) = \bigl( e_{1}^{(i)}, e_{1}^{(i)} \bigr) = (i^{2}-\sqrt{5})^{2} . \]
Then $e_{0}^{(i)}$ and $e_{1}^{(i)}$ can be seen to define disjoint hyperplanes in $\Hyp^{n}$ by \eqref{eqn:disjointnesscondition}: note that the inequality is strict.

\par For our choice of $K$ there is only one nontrivial Galois automorphism $\sigma\colon a+b\sqrt{5}\mapsto a-b\sqrt{5}$, so we also compute
\[ \bigl( e_{0}^{(i)} , e_{1}^{(i)} \bigr)^{\sigma} = i^{4}-5 \quad\text{ and }\quad \bigl( e_{0}^{(i)}, e_{0}^{(i)} \bigr)^{\sigma} = \bigl( e_{1}^{(i)}, e_{1}^{(i)} \bigr)^{\sigma} = (i^{2}+\sqrt{5})^{2} .\]
The proof of Lemma \ref{lemma} gives two ideals $\p_{0}^{(i)}= \bigl( 2 \bigl( e_{0}^{(i)},e_{0}^{(i)} \bigr) \bigr)$ and $\p_{1}^{(i)} = \bigl( 2C \bigl( e_{1}^{(i)},e_{1}^{(i)} \bigr) \bigr)$, and we require that
\[ C \geq \frac{(i^{2}+\sqrt{5})}{(i^{2}-\sqrt{5})} \qquad\qquad\text{(cf.\ \eqref{eqn:Cchoice})} . \]
We also need $e_{0}^{(i)}\pm e_{1}^{(i)}$ to be nonzero modulo $\p_{1}^{(i)}$.
That is,
\[ (-2i, 0, \ldots, 0, 2i^{2}) \quad\text{ and }\quad (2i, 0, \ldots, 0, -2\sqrt{5}) \]
are not zero modulo $\p_{1}^{(i)}$.
Since $\p_{1}^{(i)}= \bigl(2C (i^{2}-\sqrt{5})^{2} \bigr)$, this holds automatically.
Observe that if $i$ is large, then $C=2$ is sufficient.
Note also that since $\p_{0}^{(i)}$ divides $\p_{1}^{(i)}$, we need only consider $\p_{1}^{(i)}$ and can take $\Gamma'_{i}=\Gamma(\p_{1}^{(i)})$.
(Actually by Remark \ref{rem:strictness} we needn't make this justification but it is included here for completeness of exposition.)

\par Note that the proof of Theorem \ref{thm:existence} requires Lemma 3.1 to be applied a second time, with $e_{0}^{(i)}$ and $e_{1}^{(i)}$ interchanged.
However, since both vectors are of the same length, the ideal $\mathcal{Q}_{1}^{(i)}=\bigl( 4 (e_{0}^{(i)},e_{0}^{(i)} ) \bigr)$ is equal to $\p_{1}^{(i)}$ anyway, and so we can effectively ignore this step.

\par Now from hyperbolic geometry, we have \cite[Theorem~3.2.8]{Ratcliffe:Foundations}
\begin{equation}\label{eqn:distance}
  \cosh\rho\bigl(H_{0}^{(i)},H_{1}^{(i)}\bigr) = \frac{\bigl| \bigl(e_{0}^{(i)},e_{1}^{(i)}\bigr)\bigr|}{\bigl\|e_{0}^{(i)} \bigr\| \bigl\|e_{1}^{(i)}\bigr\|} = \frac{i^{2}+\sqrt{5}}{i^{2} - \sqrt{5}} 
\end{equation}
where $\rho$ is the distance between the hyperplanes $H_{0}^{(i)}$ and $H_{1}^{(i)}$ defined by $e_{0}^{(i)}$ and $e_{1}^{(i)}$ respectively.
We see that $\rho\to 0$ as $i\to\infty$.

\par In the manifold $M_{i}$ obtained by the inbreeding construction, we have $\epsilon_{i}=\Syst(M_{i}) = 2\rho_{i}$ where $\rho_{i}= \rho\bigl(H_{0}^{(i)},H_{1}^{(i)}\bigr)$.
Now, by \eqref{eqn:distance}, 
\[ \epsilon_{i} = 2\cosh^{-1} \left( \frac{i^2 + \sqrt{5}}{i^{2} - \sqrt{5}}\right) .\]
Using a Taylor expansion we deduce that for some constant $\delta>0$ we have 
\[ \epsilon_{i} \sim \frac{\delta}{i} \quad\text{ for large $i$}. \]
Writing $\p_{1}^{(i)}=(\beta)$ with $\beta=4(i^{2}-\sqrt{5})^{2}$, we have $|N(\p_{1}^{(i)})| = 16(i^{2} - \sqrt{5})^{2} (i^{2} + \sqrt{5})^{2} \sim B(\delta/\epsilon_i)^8$ for some constant $B>0$.

\par Now for a given $\p^{(i)}$, note that $|N(\p^{(i)})|$ is the number of elements in the residue class ring $\Int_{K}/\p^{(i)}$ \cite[Ch.~I, Sect.~7]{Lang:ANT}, so $|\Gamma:\Gamma'_i|\leq |N(\p^{(i)})|^{(n+1)^2}$ since $|\Gamma:\Gamma'_{i}|$ is the order of a matrix group over $\Int_{K}/\p^{(i)}$.
Thus for some positive constant $D$, 
\[ \Vol(\Hyp^{n}/\Gamma'_i) = \Vol(\Hyp^{n}/\Gamma)\cdot|\Gamma:\Gamma_{i}'| \leq D \bigl( B (\delta/\epsilon_{i})^{8}\bigr)^{(n+1)^2}, \]
which is a polynomial in $1/\Syst(M_i)$ of degree $8(n+1)^2$.
\end{proof}

\par Other related constructions of extremal sequences of manifolds $\{M_i\}$ are given in detail in the second author's forthcoming Ph.D.\ thesis \cite{Thomson:PhDthesis}.

\section{Remarks}\label{sec:5}

\subsection{Arithmeticity} If $\epsilon$ in Theorem~\ref{thm:existence} is less than some $\epsilon_0 > 0$, which depends only on the degree $d$ of the field $K$ in the construction, and on the dimension $n$, then the manifolds $M$ produced by the theorem are \emph{non-arithmetic}.
This can be seen as follows:
Assume that the manifold $M$ is arithmetic.
The fundamental group $\pi_{1}(M')$ (of the compact manifold $M'$ with boundary, of which $M$ is a double) injects into $\pi_{1}(M)$ and both are Zariski dense in $\PO(n,1)^{\circ}$.
This is shown in \cite[Sects 0.2 and 1.7]{GPS:nonarith}.
By the commensurability criterion \cite[Sect.~1.6]{GPS:nonarith} we conclude that $\Gamma_{M}=\pi_{1}(M)$ is commensurable with $\PO_{f}(\Int_{K})$.
Now we can follow a known argument relating the lengths of geodesics of an arithmetic manifold $M$ to eigenvalues of integral matrices (see \cite[Sect.~10]{Gelander}).
This implies in our case that $\Syst(M)\geq C_{n,d}$, since the integral polynomials which arise have their degree bounded above by $d(n+1)$.
Hence if $\epsilon < C_{n,d}$, then $M$ has to be non-arithmetic.
\par A conjecture of Lehmer from number theory claims that there exists a constant $m>1$ such that the Mahler measure $M(P)$ of any non-cyclotomic integral monic polynomial $P$ satisfies $M(P)\geq m$.
(Recall that the Mahler measure of an integral monic polynomial $P$ of degree $d$ with roots $\theta_{1},\theta_{2},\ldots,\theta_{d}$ is defined by $M(P)=\prod_{i=1}^{d} \max(1,|\theta_{i}|)$.)
Our argument shows that if Lehmer's conjecture is true then $\epsilon_{0}$ in the statement above is an absolute constant which does not depend on $K$ or $n$.
We refer to \cite[p.~322]{Margulis:book} and \cite[Sect.~10]{Gelander} for a related discussion addressing arithmetic manifolds.

We have shown that our method provides a new construction of non-arithmetic hyperbolic $n$-manifolds for every dimension $n$. It has some similarities with the interbreeding construction of Gromov and Piatetski-Shapiro \cite{GPS:nonarith} but at the same time is different from the former. Following Agol \cite{Agol:systoles} we call it an \emph{inbreeding construction}.

\subsection{Commensurability} If $\epsilon \to 0$ then at most finitely many of the manifolds $M$ provided by Theorem~\ref{thm:existence} will be commensurable to each other. Indeed, assume that we have an infinite sequence of commensurable non-arithmetic manifolds $M_1 = \Lambda_1\backslash\Hyp^n$, $M_2 = \Lambda_2\backslash\Hyp^n$, \ldots,\ such that $\Syst(M_i) = \epsilon_i \to 0$ when $i\to\infty$. By Margulis' Theorem \cite[Theorem~1, p.~2]{Margulis:book}, the commensurability group $\Gamma$ of $\Lambda_i$ will be a lattice in $\isom$ and hence $\Gamma \backslash \Hyp^{n}$ is a compact hyperbolic $n$-orbifold.
We have
\[\Gamma \supset \Lambda_1, \Lambda_2, \ldots, \]
so the orbifold $\Gamma\backslash\Hyp^n$ has systoles of arbitrarily short length, which is impossible. Note that this argument works for any non-arithmetic manifolds with short geodesics, not only those provided by our theorem.

\subsection{Non-compact manifolds}
The notion of systole being defined as the length of a shortest closed geodesic in a manifold $M$ generalises to non-compact finite volume hyperbolic $n$-manifolds.
In terms of the lattice $\Lambda\subset\Isom(\Hyp^{n})$ uniformising $M$, the closed geodesics in $M$ correspond to the hyperbolic elements of $\Lambda$, while $\Lambda$ also contains parabolics which have zero displacement and correspond to the cusps.
With these observations at hand the results of this paper can be generalised to the finite volume non-compact hyperbolic $n$-manifolds.
The proofs are entirely similar and we omit them.

\subsection{Some applications}
Our non-arithmetic manifolds $M$ contain properly embedded separating totally geodesic hypersurfaces, and hence the fundamental group $\pi_{1}(M)$ has the structure of a free product with amalgamated subgroup similar to the fundamental groups of the Gromov-Piatetski-Shapiro manifolds \cite{GPS:nonarith}.
This enables one to use our manifolds for the construction of Belolipetsky and Lubotzky \cite{BL:isometries}, which proves that every finite group can be realised as the full isometry group of some compact hyperbolic $n$-manifold.
Another immediate application is to the construction of some new non-coherent lattices in $\Isom(\Hyp^{n})$, which can be achieved by applying the argument of Kapovich-Potyagailo-Vinberg \cite[Sect.~4]{KPV:noncoherence} to the lattices provided by our construction.
The details of this application are explained in \cite{Thomson:PhDthesis}.

\subsection{Other locally symmetric spaces and a conjecture of Lehmer}
It is natural to ask what can be said about systoles of other locally symmetric manifolds.
This question pertains to the quotients of symmetric spaces $X=H/K$ by torsion-free lattices, where $H$ is now a semisimple Lie group and $K$ its maximal compact subgroup.

\par If all lattices in $H$ are arithmetic and Lehmer's conjecture (or its weaker version by Margulis \cite[Ch.~IX, Sect.~4.21]{Margulis:book}) holds, then the systoles of the $X$-locally symmetric manifolds will be bounded below by a constant which depends only on $X$ (cf.~\cite[Sect.~10]{Gelander}).
The arithmeticity of lattices is known for groups of real rank at least 2 by Margulis \cite[Theorem~1, p.~2]{Margulis:book}, and for
$H=\Sp(n,1)$ or $\mathrm{F}^{-20}_4$ by Corlette \cite{Corlette}. Hence the only case for which one may hope to have a version of our result is when $H=\PU(n,1)$ and $X$ is complex hyperbolic $n$-space.

There is also a reverse connection between Lehmer's conjecture and systoles of arithmetic locally symmetric manifolds. This is explained in detail for the low dimensional hyperbolic manifolds in \cite[Ch.~12.3]{MR:book} and \cite{B:report}.

\providecommand{\bysame}{\leavevmode\hbox to3em{\hrulefill}\thinspace}
\providecommand{\MR}{\relax\ifhmode\unskip\space\fi MR }
\providecommand{\MRhref}[2]{%
  \href{http://www.ams.org/mathscinet-getitem?mr=#1}{#2}
}
\providecommand{\href}[2]{#2}


\end{document}